%% file: main_v13-short-arxiv.tex
\documentclass[twocolumn]{autart2}   

\usepackage{amssymb}
\usepackage{amsmath}
\usepackage{amsfonts}
\usepackage{indentfirst}
\usepackage{lipsum}[1]
\usepackage[latin1] {inputenc}
\usepackage{enumerate}
\usepackage{mathrsfs}
\usepackage{tabulary}
\usepackage{cite}
\usepackage{wrapfig}
\usepackage{psfrag}
\usepackage{graphicx}          
\usepackage{color}      
\usepackage{epsfig} 
\usepackage{times} 
\usepackage{bbold}
\usepackage{dsfont}
\usepackage{yhmath}
\input{psfig.tex}
\input{latinte.tex}

\usepackage[prependcaption,colorinlistoftodos]{todonotes}

\def\be{\begin{equation}}
\def\ee{\end{equation}}
\def\ba{\begin{array}}
\def\ea{\end{array}}
\def\eqa{\begin{eqnarray}}
\def\eqe{\end{eqnarray}}

\begin{document}

\begin{frontmatter}

\title{Distributed averaging integral Nash equilibrium seeking on networks}

%\sgcomment{If we remove ``algorithms'' from the title, we save 1 line}

\author[claudio]{Claudio De Persis}\ead{c.de.persis@tug.nl}, \ 
\author[sergio]{Sergio Grammatico}\ead{s.grammatico@tudelft.nl}
\thanks{This work was partially supported by NWO under research projects URSES-ENBARK (408.13.037), OMEGA (613.001.702) and P2P-TALES (647.003.003), and by the ERC under research project COSMOS (802348).}

\address[claudio]{ENTEG and J.C.~Willems Center for Systems \& Control, University of Groningen, The Netherlands}
\address[sergio]{Delft Center for Systems and Control, TU Delft, The Netherlands}

\begin{abstract}
%It has been recently shown  that 
Continuous-time gradient-based Nash equilibrium seeking algorithms enjoy a passivity property under a suitable monotonicity assumption. This feature has been exploited to design distributed algorithms that converge to Nash equilibria and use local information only.  We further exploit the passivity property to interconnect the algorithms with distributed averaging integral controllers that tune on-line the weights of the communication graph. {The main advantage is to guarantee convergence to a Nash equilibrium without requiring a strong coupling condition on the algebraic connectivity of the communication graph over which the players exchange information, nor a global high-gain.}
%An advantage of these controllers is to guarantee convergence to a Nash equilibrium without requiring a strong coupling condition on the algebraic connectivity of the communication graph over which the players exchange information, nor a global high-gain. 
%Due to the presence of constraints, the stability proof relies on an invariance principle for projected dynamical systems. 
\end{abstract}

\end{frontmatter}

\section{Introduction}
In the literature on optimization and control, continuous-time algorithms that solve optimization problems are gaining increasing attention for several reasons
\cite{cherukuri.scl16,goebel.scl17}. First, the control-theoretic properties 
%of these continuous-time dynamics 
are more easily unveiled in a continuous-time setting and this permits to  naturally establish connections with control methodologies, eventually leading to technical advances in the analysis and design of the optimization algorithms themselves. Second, a way to achieve optimal performance in control systems is to interconnect the physical process with optimization algorithms and study the stability and optimality of the 
%entire 
{resulting} 
closed-loop system \cite{stegink2017frequency}. 
%\sgcomment{I would remove/rephrase this paragraph below}
%\red{Since in many important cases, the model of the system to control is in continuous-time, the interconnection of the process and the control algorithm is more natural if the latter is described in continuous-time as well.}
In the case of complex network systems, such as power and social networks,  multiple {agents or} players make decisions to optimize their own objective functions. 
{This leads to a game-theoretic setup, where}
%In this game-theoretic context, 
algorithms are commonly designed with the purpose of converging to Nash equilibria, possibly using limited or local information. The local nature of the information is defined by the topology of the network over which the game takes place. With similar motivations as before, more attention is {currently} being paid to Nash equilibrium seeking algorithms 
%implemented 
in continuous-time. 

%\cmargin{I have a poor knowledge of the literature, thus this review is rather shallow}
{\em Related literature:} Nash-equilibrium seeking algorithms have a long history and solutions in continuous-time were proposed in classical early work on game-theoretic problems \cite{rosen}. Game-theoretic problems have also attracted the interest of the control community already decades ago \cite{li1987distributed}.
%already focused on continuous-time solutions and considered generalized notions of Nash equilibria arising from the presence of constraints coupling the decisions of the different players.  
%More recently, much work has been devoted to the design of  algorithms for game equilibria computation, both in discrete-time and continuous time
%\cite{marden-et-al-tac09, cortes-martinez-tac15}. 
%N. Li and J. R. Marden, ?Designing games for distributed optimization,?
%IEEE Journal of Selected Topics in Signal Processing, vol. 7, no. 2, pp.
%230?242, 2013
%J. Marden, G. Arslan, and J. Shamma, ?Joint strategy fictitious play with
%inertia for potential games,? IEEE Trans. on Automatic Control, vol. 54,
%no. 2, pp. 208?220, 2009.
%B. Gharesifard and J. Cortes, ?Distributed convergence to Nash equilibria
%in two-network zero-sum games,? Automatica, vol. 49, no. 6, pp.
%1683 ? 1692, 2013.
%A. Cortes and S. Mart ? ??nez, ?Self-triggered best-response dynamics for
%continuous games,? IEEE Trans. on Automatic Control, vol. 60, no. 4,
%pp. 1115?1120, 2015.
%J. Koshal, A. Nedic, and U. Shanbhag, ?Distributed algorithms for ?
%aggregative games on graphs,? Operations Research, vol. 64, no. 3, pp.
%680?704, 2016.
%F. Salehisadaghiani and L. Pavel, ?Distributed Nash equilibrium seeking:
%A gossip-based algorithm,? Automatica, vol. 72, pp. 209?216, 2016.
Extensions to the case of games with coupling constraints where generalized Nash equilibria are of interest, has ben  the subject of investigation since \cite{rosen}, 
with a wide variety of results available\cite{facchinei.book,pavel2007extension,yin2011nash,kulkarni2012variational,arslan2015games,belgioioso:grammatico:17}. 
%and the investigation of dynamical aspects of game-theoretic problems have a long tradition also in the systems and control literature, e.g.~\cite{basar.olsder}. 
For general $N$-player games, the implementation of  algorithms for Nash computation require each player to access information regarding their own objective functions as well as the decisions taken by the other players in the game, information which might not be available. A way to remedy this lack of information is provided by model-free methods inspired by extremum seeking algorithms, see e.g.~\cite{frihauf-extr-seeking}, \cite{stankovic2012distributed} and references therein.  Other solutions have been proposed to solve game equilibria in a distributed fashion \cite{li-marden-STSP13, gharesifard-cortes-aut13, salehisadaghiani2016distributed,
koshal2016distributed}.  
%Aggregative games 
The approach to Nash equilibrium distributed computation which is of major interest for this paper is the one suggested by \cite{pavel}, which reconstructs the non-local information concerning the other players in the game based on a communication graph where each player communicates only with its neighbors. 
{Remarkably}, from a methodological point of view, the paper \cite{pavel} has pointed out the passivity property of Nash equilibrium seeking algorithms and used this property to design and analyze a consensus-based algorithm that uses local information only. 

{\it Paper contribution:} In this paper, we further exploit the passivity property revealed in \cite{pavel} to enrich the features of Nash equilibrium seeking algorithms. We propose a new passivity-based algorithm that allows us to relax the requirements on the knowledge of the algebraic connectivity of the graph or the use of a high-gain controller, which were the solutions  proposed in \cite{pavel}. Even though the high-gain controller of \cite{pavel} is devised with the purpose of avoiding an {\em a priori} knowledge of the algebraic graph connectivity, it still uses a global parameter ({its high-enough gain}), which might be difficult to {estimate or} implement in a network system. 
%Furthermore, compared with \cite{pavel}, we consider the presence of coupling constraints.  
%, in addition to having some well-known drawbacks, such as fast transient oscillations and lack of robustness to noise. 
%Furthermore, in \cite{pavel}, the high-gain controller was only proposed for the case of games without constraints, thus excluding an important class of games. 
Motivated by \cite{kim.depersis}, the control algorithm we propose tunes on-line the weights of the graph over which the players exchange the coordinating information via a distributed integral control action, thus enabling the convergence to Nash equilibria without assuming any information on the graph nor resorting to global parameters. Distributed averaging integral algorithms  for Nash equilibria computation of the kind studied in this paper are new to the best of our knowledge. In \cite{de-persis-grammatico-ecc18}, \cite{de-persis-grammatico-note18} different integral control laws are proposed to solve aggregative games in continuous-time.  
%Solutions are provided for games where the decisions for each player are constrained to belong to a compact set
%, as well as for the case where all the players decision variables must fulfil a set of  of affine coupling constraints. In the latter case, the proposed algorithm converges to a set of generalized Nash equilibria. 
The results are discussed in a crescendo of complexity, from {unconstrained games} (\S \ref{sec:dai.nash}), to {games with local constraints (\S \ref{sec:proj.dai.nash})}. 
%, to finish with the case of games with additional coupling constraints.  
%This allows us 
{We prefer in fact} to present the main idea in a more accessible context {first}, and then extend it to a more technical setting. 

%After recalling some technical preliminaries used throughout the paper in Section \ref{sec:prelim}, we introduce our distributed averaging controller in Section \ref{sec:dai.nash}, which is then revised to take into account local constraint sets 
%%and coupling constrains 
%in Section \ref{sec:proj.dai.nash}.  
%%and \ref{sec:coupl.constr}, respectively.

\textit{Basic notation:} 
$\R$ denotes the set of real numbers. 
%and $\overline{\R} := \R \cup \{\infty\}$ the set of extended real numbers. 
% WE NEVER USE EXTENDED REALS
$I$ denotes the identity matrix. $\bs{0}$ ($\bs{1}$) denotes a matrix/vector with all elements equal to $0$ ($1$); to improve clarity, we may add the dimension of these matrices/vectors as subscript. 
$A \otimes B$ denotes the Kronecker product between matrices $A$ and $B$. $\left\| A \right\|$ denotes the maximum singular value of matrix $A$.
Given $N$ vectors $x_1, \ldots, x_N$, we define 
%$\boldsymbol{x} := {\rm col}\left(x_1,\ldots,x_N\right) = \left[ x_1^\top \cdots x_N^\top \right]^\top$, 
${x} := {\rm col}\left(x_1,\ldots,x_N\right) = \left[ x_1^\top \cdots x_N^\top \right]^\top$, 
and, for each $i$, 
%$\boldsymbol{x}_{-i} := {\rm col}\left(x_1,\ldots , x_{i-1}, x_{i+1},\ldots,x_N\right)$.
${x}_{-i} := {\rm col}\left(x_1,\ldots , x_{i-1}, x_{i+1},\ldots,x_N\right)$.
Similarly, given $N$ sets $\Omega_1, \ldots, \Omega_N$, we define $\bs{\Omega} := \Omega_1 \times \ldots \times \Omega_N$ and, for each $i$, $\bs{\Omega}_{-i}:=\Omega_1\times\ldots\times \Omega_{i-1}\times \Omega_{i+1}\times\ldots \times \Omega_N$. For a differentiable function $v:\R^n \rightarrow \R$, 
we denote the vector of partial derivatives as $\nabla v(x) := {\rm col}\left( \tfrac{\partial v(x)}{\partial x_1}, \ldots, \tfrac{\partial v(x)}{\partial x_N}  \right) \in \R^n$. \\
\textit{Operator-theoretic definitions: }
The mapping $\iota_{\mc{S}}:\R^n \rightarrow \{ 0, \, \infty \}$ denotes the indicator function for the set $\mc{S} \subseteq \R^n$, i.e., $\iota_{\mc{S}}(x) = 0$ if $x \in \mathcal{S}$, $\infty$ otherwise. The set-valued mapping $N_{\mc{S}}: \R^n \rightrightarrows \R^n$ denotes the normal cone operator for the set $\mc{S} \subseteq \R^n$, i.e., 
$N_{\mc{S}}(x) = \emptyset$ if $x \notin \mc{S}$, $\left\{ v \in \R^n \mid \sup_{z \in \mc{S}} \, v^\top (z-x) \leq 0  \right\}$ otherwise. The set-valued mapping $T_{\mc{S}}: \R^n \rightrightarrows \R^n$ denotes the tangent cone operator for the set $\mc{S} \subseteq \R^n$.
%For a function $f: \R^n \rightarrow \overline{\R}$, $\dom(f) := \{x \in \R^n \mid f(x) < \infty\}$; $\partial f: \dom(f) \rightrightarrows {\R}^n$ denotes its subdifferential set-valued mapping, defined as $\partial f(x) := \{ v \in \R^n \mid f(z) \geq f(x) + v^\top (z-x)  \textup{ for all } z \in \textup{dom}(f) \}$; if $f$ is differentiable, then $\partial f(x) = \left\{ \nabla f(x)\right\}$, for all $x \in \R^n$.
A mapping $F:\R^n \rightarrow \R^n$ is $\ell$-Lipschitz continuous if, for all $x, y \in \dom(F)$, $\left\| F(x)-F(y)\right\| \leq \ell \left\| x-y\right\|$. A mapping $F:\R^n \rightarrow \R^n$ is (strictly) monotone if, for all $x, y \in \dom(F)$, $\left( x - y\right)^\top \left( F(x)-F(y)\right) (>) \geq 0$; $F$ is $\epsilon$-strongly monotone, with $\epsilon>0$, if, for all $x, y \in \dom(F)$, $\left( x - y\right)^\top \left( F(x)-F(y)\right) \geq \epsilon \left\| x-y\right\|^2$.
Given a closed convex set $\mathcal{C} \subseteq \R^n$ and a single-valued mapping $F: \mathcal{C} \rightarrow \R^n$, the variational inequality problem VI$(\mathcal{C},F)$, is the problem to find $x^* \in \mathcal{C}$ such that 
$(y-x^*)^\top  \, F(x^*) \geq 0 \ \textup{  for all } y\in \mathcal{C}.$
The mapping $\proj_{\Omega}$ denotes the projection operator for the set $\Omega \subseteq \R^n$, i.e., $\proj_{\Omega}(x) := {\rm argmin}_{ y \in \Omega } \left\| y-x\right\| $; $\Pi_{\Omega}(x,v)$ denotes the projection of the vector $v$ onto the tangent cone of $\Omega$ at $x$, i.e., $\Pi_{\Omega}(x,v) := \proj_{ T(x) }(v)$.

\section{Mathematical {background}}\label{sec:prelim}

In this section, we recall a few known results to set the ground for the analysis in the rest of the paper. 
%\sgcomment{I would emphasize that this section is a background section in the sense that the results given are known, although recent.}

%\subsection{Game-theoretic background}
%\subsection{{Nash equilibrium seeking via static consensus}}
\subsection{{Nash equilibrium problem}}
A game, denoted by $G(\mathcal{I}, \{ J_i \}_i, \{ \Omega_i \}_i)$ consists of $N$ players indexed by the set $\mathcal{I}:=\{1,2,\ldots, N\}$, 
where each agent $i$ can decide on a strategy vector $x_i\in \Omega_i \subseteq \mathbb{R}^{n_i}$, with the aim to minimize its the cost function $J_i (x_i, {x}_{-i})$. Let us define $M := \sum_{i=1}^{N} n_i \geq N$.

%\cmargin{the notation $\bs{x}^*$ to denote ${\rm col}(x_1^*,...,x_N^*)$, with  $\bs{x}^*\in \mathbb{R}^M$, gets confused with the symbol $\bs{x}$, which we use to denote 
% ${\rm col}(\bs{x}^1,...,\bs{x}_N)\in \mathbb{R}^{NM}$.}
\begin{definition}[Nash equilibrium]
A vector of strategies ${x}^* \in \bs{\Omega} \subseteq \R^M$ is a Nash equilibrium if, for all $i \in \mathcal{I}$,  
\[
J_i(x_i^*, x_{-i}^*) \le \inf_{ x_i \in \Omega_i } J_i(x_i, x_{-i}^*).
\vspace{-0.8cm}
\]
\hfill $\square$
\end{definition}

\begin{assumption}[Convexity]
\label{strict.convexity}{\rm \cite[Assumption 2]{pavel}}\\
(i) For each $i\in \mathcal{I}$, $\Omega_i =\mathbb{R}^{n_i}$, and $J_i(x_i, x_{-i})$ is $\mathcal{C}^2$, strictly convex and radially unbounded in $x_i$ for every ${x}_{-i} \in \bs{\Omega}_{-i}$.  \\
(ii) For each $i\in \mathcal{I}$, $\Omega_i \subset \mathbb{R}^{n_i}$ is a nonempty, convex and compact set, $J_i(x_i, x_{-i})$ is $\mathcal{C}^1$ and convex in $x_i$ for every ${x}_{-i}\in \bs{\Omega}_{-i}$. 
\hfill $\square$
\end{assumption}

It follows from \cite[Cor. 4.2]{basar.olsder} that under Assumption \ref{strict.convexity}(i), a Nash equilibrium $x^*$ exists, and it satisfies
$\frac{\partial J_i}{\partial x_i}(x_i^*, {x}_{-i}^*)=0$, for all $i \in \mathcal{I}$.
In fact, for the existence of a Nash equilibrium, continuity of  $J_i$ is sufficient \cite[Th. 1]{rosen}. 

%\sgmargin{x no bold (red)?}

Let us introduce the pseudo-gradient dynamics:
\be\label{pseudo-gradient-dynamics}
\dot{{x}} = 
\begin{bmatrix}
\dot x_1\\
\vdots \\
\dot x_N
\end{bmatrix}
=
-
\begin{bmatrix}
\dst\frac{\partial J_1}{\partial x_1}(x_1, {x}_{-1})\\[3mm]
\vdots \\
\dst\frac{\partial J_N}{\partial x_N}(x_N, {x}_{-N})\\[3mm]
\end{bmatrix}
=: -F({x}).
\ee

\begin{assumption}[Strictly monotone pseudo-gradient]
\label{strict.monotonicity}
The pseudo-gradient mapping $F$ in \eqref{pseudo-gradient-dynamics} is strictly monotone.
\hfill $\square$
\end{assumption}

Due to Assumption \ref{strict.monotonicity}, it follows from \cite[Th. 3]{scutari} that there exists a unique Nash equilibrium for the game $G(\mathcal{I}, \{ J_i \}_i, \{\Omega_i\}_i)$. Moreover, due to Assumptions \ref{strict.convexity} and  \ref{strict.monotonicity}, the equilibrium of the dynamics in \eqref{pseudo-gradient-dynamics} is the unique Nash equilibrium and is globally asymptotically stable.

\subsection{{Nash equilibrium problem with partial information: Distributed seeking via static consensus}}

Whenever an agent $i$ has no access to the full information ${x}_{-i}$, the authors in \cite{pavel} propose to augment the pseudo-gradient dynamics. Specifically, each agent shall estimate the states of all the other agents, i.e., implement the following dynamics:
\be\label{eq:ext.pseudo.dynamics}
\ba{rcl}
\dot x_j^{i} &=& u_j^i \quad \forall j \neq i \\ [3mm]
\dot x_{i} &=& -\dst\frac{\partial J_i}{\partial x_i} \left(x_i, \bs{x}_{-i}^i\right) + u_{i}^i \\
\ea
\ee
where {$\bs{x}^i:={\rm col}(x_1^{i}, \ldots, x_{i-1}^{i}, x_i,  x_{i+1}^{i}, \ldots, x_N^{i})$ is the state vector of agent $i$,} 
and $u_j^i\in \mathbb{R}^{{n_j}}$, for all $j \in \mathcal{I}$, are the control inputs for agent $i$ to be designed. Note that in \eqref{eq:ext.pseudo.dynamics}, 
 each agent $i$ relies on the estimated strategies $\bs{x}_{-i}^i$, which is available locally, not on the true strategies $\bs{x}_{-i}$. This approach comes at the expenses that each agent has a local copy of the state of each other agent in the game. Whenever confusion does not arise, {let us use} both $x_{i}^{i}$ and $x_{i}$ for the same variable. 

The dynamics for agent $i$ in compact form read as
\[
\dot{\bs{x}}^i = - \mathcal{R}_i^\top \, \dst\frac{\partial J_i}{\partial x_i} \left(\bs{x}^i\right) + \bs{u}^i,
\] 
where $\mathcal{R}_i \in \R^{ n_i \times M }$ is the matrix
\[
\mathcal{R}_i :=
\begin{bmatrix}
\bs{0}_{n_i \times {n_1}} {\cdots} \, \bs{0}_{n_i \times n_{i-1}}\, I_{{n_i\times n_i}}\, \bs{0}_{n_i \times n_{i+1}} {\cdots} \, \bs{0}_{n_i \times n_N}
\end{bmatrix}.
\]

To induce the local variables $x^i_j$, $j\ne i$, to converge towards the true values $x_j$, a  consensus protocol is used, i.e., 
\begin{equation} \label{eq:u-static-consensus}
\bs{u}^i = \textstyle \sum_{j=1}^N a_{i,j} \left( \bs{x}^j - \bs{x}^i \right),
\end{equation}
where $a_{i,j}$ are the entries of the adjacency matrix of a graph over which the agents exchange information. 

%\cmargin{Is ``undirected" needed?}
%\cmargin{posticiperei questo rilassamento}
\begin{assumption}\label{asspt:graph}
The information exchange graph is undirected and connected.  
\hfill $\square$
\end{assumption}

{Now, if we define}
$\boldsymbol{x}= \col\left(\bs{x}^1, \ldots, \bs{x}^N \right),$

\begin{equation}
\label{eq:extended.pseudo-gradient}
\boldsymbol{F}(\boldsymbol{x})= \col\left(\dst\frac{\partial J_1}{\partial x_1} (\bs{x}^1), \ldots, \dst\frac{\partial J_N}{\partial x_N} (\bs{x}^N) \right)
\end{equation}

and 
$
\mathcal{R}^\top= {\rm diag}\left( \mathcal{R}_1^\top, \ldots, \mathcal{R}_N^\top\right).
$, then, the closed-loop system in compact form reads as 
\be\label{eq:cl-loop-compact}
\dot{\boldsymbol{x}}= -\mathcal{R}^\top \boldsymbol{F}(\boldsymbol{x}) - (L \otimes I_{M}) 
\boldsymbol{x},
\ee

where $L$ is the Laplacian matrix associated with the information-exchange graph.

It follows from \cite[Lemma 2]{pavel} that an equilibrium $\bar{\boldsymbol{x}}$ of \eqref{eq:cl-loop-compact} satisfies $\bar{\boldsymbol{x}}=\bs{1}_N \otimes {x}^*$, with ${x}^*$ {being} a Nash equilibrium.

\subsection{Stability analysis}

For the stability analysis of the dynamics in \eqref{eq:cl-loop-compact}, we consider the following assumptions.

\begin{assumption}[Strongly monotone pseudo-gradient]
\label{ass:strong.monotonicity}
The pseudo-gradient mapping $F$ in \eqref{pseudo-gradient-dynamics} is $\mu$-strongly monotone, with $\mu>0$.
\hfill $\square$
\end{assumption}

%\cite[Assumption 4(ii)]{pavel}

\begin{assumption}[Lipschitz pseudo-gradient]
\label{asspt:lipschitz.extended}
The pseudo-gradient mapping $F$ in \eqref{pseudo-gradient-dynamics} is $\ell_{F}$-Lipschitz continuous, with $\ell_{F}>0$; the extended pseudo-gradient mapping $\bs{F}$ in \eqref{eq:extended.pseudo-gradient} is $\ell_{\bs{F}}$-Lipschitz continuous, with $\ell_{\textbf{F}}>0$. Define $\ell := \max\{ \ell_F, \ell_{\bs{F}} \}$.
\hfill $\square$
\end{assumption}

The following result due to \cite{pavel} establishes global asymptotic convergence of the dynamics in \eqref{eq:cl-loop-compact} to the unique Nash equilibrium, under a strong coupling {condition.}

\begin{lem} \label{th:pavel} 
{\rm\cite[Th. 2]{pavel}}
Let Assumptions \ref{strict.convexity}(i), \ref{asspt:graph}, \ref{ass:strong.monotonicity}, \ref{asspt:lipschitz.extended} hold. 
If
\be\label{strong.coupling}
\lambda_2(L) > \ell + {\ell^2}/{\mu},
\ee
then, for any initial condition, the solution to \eqref{eq:cl-loop-compact} is bounded and converges {exponentially} to the {unique} Nash equilibrium, i.e., $\lim_{t \to \infty} \boldsymbol{x}(t)= \bs{1}_N \otimes {x}^*$.
\hfill $\square$
\end{lem}

\section{Distributed averaging integral Nash equilibrium seeking}\label{sec:dai.nash}

In this section, we propose a novel Nash equilibrium seeking algorithm on networks. 
%under Assumption \ref{strict.convexity}(i). 
Instead of the static consensus coupling $\bs{u} = - (L \otimes I_{M}) \bs{x}$ {in \eqref{eq:u-static-consensus}}, we propose an integral loop that tunes the control gains to compensate for the possible lack of strong coupling, i.e., without requiring condition \eqref{strong.coupling} on the algebraic connectivity of the graph:

\be\label{eq:new.control.player.i}
\ba{rcl}
\dot k_i &=& \gamma_i \, \left\| \bs{\rho}^i \right\|^2 \\
\bs{\rho}^i &=& \dst \sum_{j=1}^N a_{i,j} \left( \bs{x}^j - \bs{x}^i \right) \\
\bs{u}^i &=& \dst -\sum_{j=1}^N a_{i,j}  (k_j \bs{\rho}^j - k_i \bs{\rho}^i)
\ea
\ee
where, for all $i \in \mathcal{I}$, $k_i\in \mathbb{R}$ is the state of the $i$th integrator, $\bs{u}^i, \bs{\rho}^i\in \mathbb{R}^{M}$, and  $\gamma_i>0$ is a constant parameter. 
Note that the control algorithm in \eqref{eq:new.control.player.i} requires the agents to exchange the variables $\bs{x}^i$, $k_i$, and $\bs{\rho}^i$.

In vector form, we have
\be\label{eq:new.control}\ba{rcl}
\dot{\bs{k}} &=& D(\bs{\rho})^\top (\Gamma\otimes I_{M}) \bs{\rho} \\
\bs{u} & = & -( L K L\otimes I_{M})\boldsymbol{x} 
\ea
\ee
where
\[\ba{ll}
\bs{k} := {\rm col}(k_1, \ldots, k_N),& \bs{\rho} := {\rm col}( \bs{\rho}^1, \ldots, \bs{\rho}^N), \\ 
\Gamma := {\rm diag}(\gamma_1, \ldots, \gamma_N)& K:={\rm diag}(k_1, \ldots, k_N),\\
D(\bs{\rho}) := {\rm {block.}diag}( \bs{\rho}^1, \ldots, \bs{\rho}^N),  &
\ea\]
and we used that $\bs{\rho} = -(L\otimes I_{M})\boldsymbol{x}$ and $\bs{u} = (L\otimes I_{M})(K \otimes I_{M}) \bs{\rho}$. 
Then, the resulting closed-loop system has the form
\be\label{eq:cl-loop-compact-new}
\ba{rcl}
\dot{\boldsymbol{x}} &=&  -\mathcal{R}^\top \boldsymbol{F}(\boldsymbol{x}) -(LKL\otimes I_{M})\boldsymbol{x}  \\
\dot{\bs{k}} &=& D(\bs{\rho})^\top (\Gamma\otimes I_{M}) \bs{\rho}, \quad  \bs{\rho} = -(L\otimes I_{M})\boldsymbol{x}.
\ea
\ee

Next, inspired by \cite{kim.depersis}, we show convergence of $\bs{x}(t)$ in \eqref{eq:cl-loop-compact-new} to the Nash equilibrium without the assumption of strong coupling of the information exchange graph. 

\begin{theorem}[Convergence to Nash equilibrium]
\label{convergence1}
Let Assumptions \ref{strict.convexity}(i), \ref{strict.monotonicity}, \ref{asspt:graph}, \ref{asspt:lipschitz.extended}, hold. Then, for any initial condition, the solution to system \eqref{eq:cl-loop-compact-new} is bounded and its ${\boldsymbol{x}}$-component globally asymptotically converges to the Nash equilibrium, i.e., 
$\lim_{t\to \infty} \boldsymbol{x}(t)= \bs{1}_N \otimes {x}^*$.
\hfill $\square$
\end{theorem}

\begin{pf}
See Appendix \ref{app:proof-th1}.
\hfill $\blacksquare$
\end{pf}

\section{Projected distributed averaging integral Nash equilibrium seeking}\label{sec:proj.dai.nash}
\label{sec.constrained}
In this section, we postulate Assumption \ref{strict.convexity}(ii), that is, we consider the case where $\Omega_i \subset \R^{n_i}$ for all $i$. Under Assumptions \ref{strict.convexity}(ii) and \ref{strict.monotonicity}, there exists a unique Nash equilibrium \cite{rosen, scutari,pavel}. 
Moreover, by \cite[Prop. 1.4.2]{facchinei.book}, a {vector} is the Nash equilibrium if and only if it satisfies the variational inequality VI$(\bs{\Omega},F)$. 
%\cmargin{Is the sentence in red correct?}
%\cdp{Note that compared with \cite{cherukuri.scl16}, the boundedness of $\Omega_i$ is necessary to ensure the existence of a Nash equilibrium. }

Due to the presence of the constraint sets, in \cite{pavel}, the authors consider the projected pseudo-gradient dynamics
\be\label{eq:ext.pseudo.dynamics.proj}
\ba{rcl}
\dot x_j^{i} & = & u_j^i \quad \forall j \neq i \\ [3mm]
\dot x_{i} &=& \Pi_{\Omega_i}\left( x_i \, , \, -\tfrac{\partial J_i}{\partial x_i} (x_i, \bs{x}_{-i}^i)+u_{i}^i \right)\\
\ea
\ee
that in compact form read as
\be\label{eq:ext.pseudo.dynamics.proj2}
\dot{\bs{x}}^i = \mathcal{R}_i^\top \Pi_{\Omega_i}\!\left( x_i, -\tfrac{\partial J_i}{\partial x_i} ( \bs{x}^i ) \!+\! \mathcal{R}_i \bs{u}^i \right) + (I_{M} - \mathcal{R}_i^\top \mathcal{R}_i) \bs{u}^i
\ee
where
\[
(I_{M} -
\mathcal{R}_i^\top \mathcal{R}_i) \bs{u}^i =
\begin{bmatrix}
\col\left( u_1^i, \ldots, u_{i-1}^i\right) \\
\bs{0}_{n_i} \\
\col\left( u_{i+1}^{i}, \ldots, u_{N}^i\right)
\end{bmatrix}.
\]

Similarly to \eqref{eq:cl-loop-compact}, the overall extended pseudo-gradient dynamics can be written as
\be
\label{eq:overall-proj-dynamics-compact}
\dot{\boldsymbol{x}} = \mathcal{R}^\top \Pi_{\bs{\Omega}}\left(\mathcal{R}\boldsymbol{x} \, , \, 
-\boldsymbol{F}(\boldsymbol{x}) + \mathcal{R} \bs{u} \right) + (I_{M} -\mathcal{R}^\top \mathcal{R}) \bs{u}.
\ee

We recall the following equivalence between the equilibrium of \eqref{eq:overall-proj-dynamics-compact} and the Nash equilibrium:
\begin{lem}
${x}^* \in \bs{\Omega}$ is the Nash equilibrium if and only if the pair $\left(\overline{\bs{x}}, \overline{\bs{u}}\right)=(\bs{1}_N \otimes {x}^*, \bs{0}_M)$ is the equilibrium of \eqref{eq:overall-proj-dynamics-compact}, i.e., 
\be\label{eq:eq}
\boldsymbol{0} = \mathcal{R}^\top \Pi_{\bs{\Omega}}\left(\mathcal{R}\overline{\boldsymbol{x}} \, , \, 
-\boldsymbol{F}(\overline{\boldsymbol{x}})\right).
\vspace{-0.25cm}
\ee
\hfill $\square$
\end{lem}

\begin{pf}
Equation \eqref{eq:eq} holds if and only if $\bs{x}^*$  is a solution to the VI$\left( \bs{\Omega}, F\right)$. The proof follows from \cite[Prop. 1.4.2]{facchinei.book}. 
\hfill $\blacksquare$
\end{pf}

We now rephrase a key result from \cite{pavel}.
\begin{lem}\label{lem:diss.in.pds}
{\rm \cite[Lemma 8]{pavel}}
The storage function $V(\boldsymbol{x},\boldsymbol{y})= \frac{1}{2} \|\boldsymbol{x}-\boldsymbol{y}\|^2$ satisfies, for all $\boldsymbol{x},\boldsymbol{y}$ such that $\mathcal{R} \boldsymbol{x},\mathcal{R} \boldsymbol{y} \in \bs{\Omega}$ and all $\bs{u}, \bs{v}\in \mathbb{R}^M$, 
$$
\ba{c}
%\begin{bmatrix} \dst \frac{\partial V(\bs{x},\bs{y})}{\partial \boldsymbol{x}}^\top & \dst \frac{\partial V(\bs{x},\bs{y})}{\partial \boldsymbol{y}}^\top \end{bmatrix}
\nabla V(\bs{x}, \bs{y})^\top
\left[ \begin{smallmatrix} \dot{\boldsymbol{x}} \\ \dot{\boldsymbol{y}}  \end{smallmatrix}\right] \\
\le 
-(\boldsymbol{x}-\boldsymbol{y}) \mathcal{R}^\top (\boldsymbol{F}(\boldsymbol{x})-\boldsymbol{F}(\boldsymbol{y}))
+(\boldsymbol{x}-\boldsymbol{y})^\top ( \bs{u} - \bs{v} ),
\ea
$$
where $\dot{\boldsymbol{x}}$ is the right-hand side of \eqref{eq:overall-proj-dynamics-compact} and 
$\dot{\boldsymbol{y}}$ equals 
\be\label{eq:overall-proj-dynamics-compact-x'}
\mathcal{R}^\top \Pi_{\bs{\Omega}}\left( \mathcal{R}\boldsymbol{y} \, , \, -\boldsymbol{F}(\boldsymbol{y}) +\mathcal{R} \bs{v} \right) +\left( I_{M} -\mathcal{R}^\top \mathcal{R} \right) \bs{v}.
\vspace{-0.25cm}
\ee
\hfill $\square$
\end{lem}

\medskip

We remark that in Lemma \ref{lem:diss.in.pds} the dissipation inequality is intended to hold point-wise, which dispenses us to specify the notion of solution at this stage. 
The dissipation inequality highlighted in Lemma \ref{lem:diss.in.pds} is important for our purposes because it allows us to derive a Lyapunov inequality for the feedback interconnection of the project dynamical system in \eqref{eq:overall-proj-dynamics-compact} with the distributed averaging integral control in \eqref{eq:new.control}, namely for the closed-loop system:
\be\label{eq:overall-proj-dynamics-compact-cl}
\ba{rcl}
\dot{\boldsymbol{x}} &=&  \mathcal{R}^\top \Pi_{\bs{\Omega}}(\mathcal{R}\boldsymbol{x} \, , \, -\boldsymbol{F}(\boldsymbol{x}) -\mathcal{R} (LKL\otimes I_{M})\boldsymbol{x})
\\
&& 
-(I_{M} -\mathcal{R}^\top \mathcal{R})(LKL\otimes I_{M})\boldsymbol{x},\\
\dot{\bs{k}} &=& D(\boldsymbol{\rho})^\top (\Gamma\otimes I_{M}) \bs{\rho},\quad  \bs{\rho} = -(L\otimes I_{M})\boldsymbol{x}.
\ea\ee

We can show the following inequality for the closed-loop system in \eqref{eq:overall-proj-dynamics-compact-cl}.

\begin{lem}\label{lem:diss.in.pds.cl}
Let Assumption \ref{strict.convexity}(ii) hold. Then, the Lyapunov function 
\[
W(\boldsymbol{x}, \bs{k}) = \tfrac{1}{2} \| \boldsymbol{x} - \overline{\boldsymbol{x}}\|^2  + \tfrac{1}{2} 
\left\| \bs{k} - \overline{\bs{k}} \right\|_{\Gamma^{-1}}^2,
\]
where $\overline{\boldsymbol{x}} = \bs{1}_N \otimes {x}^*$, $\overline{\bs{k}}= k^* \bs{1}_N$ and $k^*\in \mathbb{R}$ to determine, satisfies, for all $\boldsymbol{x}$ such that $\mathcal{R} \boldsymbol{x}\in \bs{\Omega}$ and all $\bs{k}\in \mathbb{R}^N$, 
\begin{equation}
\ba{rl}
\nabla W( \bs{x}, \bs{k} )
\left[ \begin{smallmatrix} \dot{\boldsymbol{x}} \\ \dot{\bs{k}} \end{smallmatrix} \right] \le & -(\boldsymbol{x}-\overline{\boldsymbol{x}}) \mathcal{R}^\top (\boldsymbol{F}(\boldsymbol{x})-\boldsymbol{F}(\overline{\boldsymbol{x}})) \\
& -(\boldsymbol{x}-\overline{\boldsymbol{x}})^\top  (L {K}^* L\otimes I_{M})(\boldsymbol{x}-\overline{\boldsymbol{x}}),
\ea
\end{equation}
where ${K}^*= k^* I_{N}$ and $\left[ \begin{smallmatrix} \dot{\boldsymbol{x}} \\ \dot{\bs{k}} \end{smallmatrix} \right] $ denotes the right-hand side of \eqref{eq:overall-proj-dynamics-compact-cl}. 
\hfill $\square$
\end{lem}

\begin{pf}
Assumption \ref{strict.convexity}(ii) implies the existence of a Nash equilibrium ${x}^*\in \bs{\Omega}$, which is equivalent to $\mathcal{R} \overline{\boldsymbol{x}}\in \bs{\Omega}$. Thus, we can apply Lemma \ref{lem:diss.in.pds} with $\overline{\boldsymbol{x}}$ and $\bs{0}$ in place of $\boldsymbol{y}$ and $\bs{v}$, respectively. 
In view of \eqref{eq:eq}, we obtain that
\[\ba{rl}
\dst \frac{\partial V( \bs{x}, \overline{\bs{x}} )}{\partial \boldsymbol{x}}^\top\dot{\boldsymbol{x}} \le & 
-(\boldsymbol{x}-\overline{\boldsymbol{x}}) \mathcal{R}^\top (\boldsymbol{F}(\boldsymbol{x})-\boldsymbol{F}(\overline{\boldsymbol{x}}))\\
 & +(\boldsymbol{x}-\overline{\boldsymbol{x}})^\top \bs{u}.
\ea
\]

Now, for ${\bs{u}}= -(L K L\otimes I_{M}) \boldsymbol{x}$, since $\tfrac{\partial V}{\partial \boldsymbol{x}}=\tfrac{\partial W}{\partial \boldsymbol{x}}$, the inequality above becomes
\begin{multline}
\label{diss.in.bold.x}
\dst \frac{\partial W}{\partial \boldsymbol{x}}^\top\dot{\boldsymbol{x}} \le 
-(\boldsymbol{x}-\overline{\boldsymbol{x}}) \mathcal{R}^\top (\boldsymbol{F}(\boldsymbol{x})-\boldsymbol{F}(\overline{\boldsymbol{x}})) \\
-(\boldsymbol{x}-\overline{\boldsymbol{x}})^\top (LKL\otimes I_{M})\boldsymbol{x}.
\end{multline}

Furthermore, as in the proof of Theorem \ref{convergence1},  
\be\label{diss.in.k}
\dst \frac{\partial W}{\partial \bs{k}}^\top\dot{\bs{k}}=\boldsymbol{x}^\top ( L (K-K^*) L \otimes I_{nN} )\boldsymbol{x}.
\ee  

Since the sum of the second addend on the right-hand side of \eqref{diss.in.bold.x} and the term on the right-hand side of \eqref{diss.in.k} is $-(\boldsymbol{x}-\overline{\boldsymbol{x}})^\top  (L {K}^* L\otimes I_{M})(\boldsymbol{x}-\overline{\boldsymbol{x}})$, the thesis follows.
\hfill $\blacksquare$
\end{pf}

Similarly to Lemma \ref{lem:diss.in.pds}, we remark that the Lyapunov inequality in Lemma \ref{lem:diss.in.pds.cl} holds point-wise. 
We now use Lemma \ref{lem:diss.in.pds.cl} and an invariance principle for projected dynamical systems to infer convergence for the closed-loop system. 

We first note that the closed-loop system \eqref{eq:overall-proj-dynamics-compact-cl} can be written as a projected dynamical system: 
\be\label{sys.overall.pds}
\left[ \begin{matrix} \dot{\boldsymbol{x}} \\ \dot{\bs{k}} \end{matrix} \right]
= \Pi_{\boldsymbol{\Xi}} \left( \left[ \begin{matrix} {\boldsymbol{x}} \\ {\bs{k}} \end{matrix} \right] \, , \, 
\boldsymbol{g}(\boldsymbol{x}, \bs{k}) \right).
\ee

Specifically, let us define 
$\bs{\Omega}^1 := \Omega^1 \times \mathbb{R}^{n} \times \ldots \times \R^n$,
$\bs{\Omega}^2 := \R^n \times \Omega^2 \times \mathbb{R}^{n} \times \ldots \times \R^n$, $\ldots \,$, 
$\bs{\Omega}^N := \R^n \times \ldots \times \R^n \times \Omega^N$, the closed convex set 
$\boldsymbol{\Xi} :=  \bs{\Omega}^1 \times \ldots \times  \bs{\Omega}^N \times \mathbb{R}^N$, for all $i \in \mathcal{I}$, the mapping 
\be\label{f^i}
f^i( \bs{x}^i, \bs{u}^i):= \mathcal{R}_i^\top \frac{\partial J_i}{\partial x_i} (x_i, \bs{x}_{-i}^i)+\bs{u}^i, 
\ee
and finally the mapping 
$$
\boldsymbol{g}(\boldsymbol{x}, \bs{k}) :=
\begin{bmatrix}
f^1\left( \bs{x}^1 \, , \,- \sum_{j=1}^N a_{1,j}  (k_j \bs{\rho}^j - k_1 \bs{\rho}^1)\right)
\\ \vdots \\ 
f^N\left( \bs{x}^N \, , \, - \sum_{j=1}^N a_{N,j}  (k_j \bs{\rho}^j - k_N \bs{\rho}^N) \right)\\
\gamma_1 \left\| \bs{\rho}^1 \right\|^2 \\
\vdots\\
\gamma_N \left\| \bs{\rho}^N \right\|^2
\end{bmatrix},
$$
where we recall that, for all $i \in \mathcal{I}$, $\bs{\rho}^i = \sum_{j=1}^{N} a_{i,j} \left( \bs{x}^j - \bs{x}^i \right) $ from equation \eqref{eq:new.control.player.i}.

The projected dynamical system in \eqref{sys.overall.pds} has a discontinuous right-hand side, and  its solutions must be intended in a Carath\`{e}odory sense. It is known \cite[Theorem 2.5]{nagurney}, \cite[Proposition 2.2]{cherukuri.scl16}, that if the vector field $\boldsymbol{g}$ is Lipschitz continuous on the closed convex set $\boldsymbol{\Xi}$, then for any initial condition $(\boldsymbol{x}_0, \bs{k}_0)\in \boldsymbol{\Xi}$, there exists a unique  Carath\`{e}odory solution to \eqref{sys.overall.pds} from $(\boldsymbol{x}_0, \bs{k}_0)$ that is defined on the entire interval $[0, \infty)$, satisfies $(\boldsymbol{x}(t), \bs{k}(t)) \in \boldsymbol{\Xi}$ for all $t \in [0, \infty)$, and is uniformly continuous with respect to the initial condition. The latter property is essential to prove the invariance of the positive limit set associated to any initial condition, as stated in the following result.

\begin{lem}\label{lem:inv.limit.set}{\rm \cite[Proof Th. 4, Part 1]{brogliato}}
If the mapping $\boldsymbol{g}$ is Lipschitz continuous on the closed convex set $\boldsymbol{\Xi}$, then, for any $\boldsymbol{y}_0\in \boldsymbol{\Xi}$, the positive limit set $\Lambda(\boldsymbol{y}_0)$ of the solution to \eqref{sys.overall.pds} starting from $\boldsymbol{y}_0$ is forward invariant.  
\hfill $\square$
\end{lem}
\smallskip

\begin{pf}
Take any point $\bs{z} \in \Lambda(\boldsymbol{y}_0)$. By definition, there exists a diverging sequence $\{t_k\}_{k\in \mathbb{N}}$
such that $\lim_{k\to \infty} \boldsymbol{y}(t_k, \boldsymbol{y}_0)= \bs{z}$. By the Lipschitz continuity, uniform continuity of the solution with respect to the initial condition holds, and therefore $\boldsymbol{y}(t, \bs{z}) =  \lim_{k \to \infty} \boldsymbol{y}(t, \boldsymbol{y}(t_k, \boldsymbol{y}_0))$. Then, one can proceed as in \cite[Proof of Th. 4, Part 1]{brogliato}. 
\hfill \blacksquare
\end{pf}

Once we have guaranteed invariance of the limit set, we establish the following invariance principle for systems with Carath\`{e}odory solutions. {The proof is a variation of the arguments in \cite{brogliato,cherukuri.scl16,bacciotti-ceragioli} adapted to our case study and is included for the sake of completeness. }

\smallskip
\begin{theorem}\label{invariance.pds}
Consider a projected dynamical system $\dot x= \Pi_K(f(x))$, where the set $K\subset \mathbb{R}^n$ is closed and convex, and the mapping $f: K \to \mathbb{R}^n$ is continuously differentiable. 
Let $\Psi\subset \mathbb{R}^n$ be a compact set such that the intersection $K\cap \Psi$ is an invariant set for $\dot x = \Pi_K(f(x))$.
Suppose that there exists a continuously differentiable function $V:\mathbb{R}^n \to \mathbb{R}$ such that 
\be\label{lasalle.ineq}
\sup_{ x\in K\cap \Psi } \dst \nabla V(x)^\top \Pi_K(f(x)) \le 0. 
\ee

Then, for any initial condition in $K\cap \Psi$, there exists a unique Carath\`{e}odory solution to $\dot x= \Pi_K(f(x))$, which remains in $K\cap \Psi$  and converges to the largest invariant set  contained in $\{ x \in K\cap \Psi\,\mid\, \nabla V(x)^\top \Pi_K(f(x))=0\}$.
\hfill $\square$
\end{theorem}

\begin{pf}
See Appendix \ref{app:th2}.
\hfill \blacksquare
\end{pf}

We are ready to show the main result of this section, the global asymptotic convergence to a Nash equilibrium of the projected dynamical system in \eqref{sys.overall.pds}. Technically, we achieve the same convergence result as in \cite[Th. 5]{pavel} without assuming a lower bound on the algebraic connectivity of the graph. 
Our proof relies on an invariance principle for projected dynamical systems, not on Barbalat's lemma as in \cite[Proof of Th. 5]{pavel}.

\begin{theorem}\label{thm:conv.constr}
Under Assumptions \ref{strict.convexity}(ii), \ref{asspt:graph}, \ref{ass:strong.monotonicity}, \ref{asspt:lipschitz.extended},  
for any initial condition 
%in $\boldsymbol{\Omega} \times \R^{N}$, 
in $\boldsymbol{\Xi}$, 
there exists a unique Carath\`{e}odory solution to system  \eqref{eq:overall-proj-dynamics-compact-cl}, which belongs to 
%$\boldsymbol{\Omega} \times \R^{N}$ 
$\boldsymbol{\Xi}$
for all $t\ge 0$, such that its $\boldsymbol{x}$-component converges asymptotically to 
the Nash equilibrium, i.e., $\lim_{t \to \infty} \boldsymbol{x}(t)= \bs{1}_N \otimes {x}^*$.
\hfill $\square$
\end{theorem}

\begin{pf}
See Appendix \ref{app:th3}.
\hfill \blacksquare
\end{pf}

\section{Discussion}

\subsection{On the proposed distributed averaging algorithm}

The proposed  algorithm in \eqref{eq:new.control.player.i} comprises  a distributed  integral action averaging the local decisions of the agents whose state variable $k$ is used as a tuning weight of the input $u$ coupling the agents' dynamics. Consequently, the proposed algorithm requires the exchange of the vector $k_i \bs{\rho}^i$ in addition to the vector $\bs{x}^i$. By following \cite{kim.depersis}, an alternative control algorithm may be given by 
\begin{equation}\label{eq:new.control.player.i.new}
\dot k_i = \gamma_i \left\| \bs{\rho}^i \right\|^2, \ \ 
\bs{\rho}^i = \dst\sum_{j=1}^N a_{i,j} ( \bs{x}^j - \bs{x}^i ), \ \ 
\bs{u}^i = k_i \bs{\rho}^i
\end{equation}
%\be\label{eq:new.control.player.i.new}
%\ba{rcl}
%\dot k_i &=& \gamma_i \left\| \bs{\rho}^i \right\|^2 \\
%\bs{\rho}^i &=& \dst\sum_{j=1}^N a_{i,j} ( \bs{x}^j - \bs{x}^i )\\
%\bs{u}^i & = & k_i \bs{\rho}^i, 
%\ea
%\ee
where the local control $\bs{u}^i$ uses the local average $\bs{\rho}^i$ only, as opposed to \eqref{eq:new.control.player.i}. This alternative approach requires the use of a different Lyapunov function and to establish additional technical results for the boundedness of the solutions. {Similarly, the results {in} \cite{kim.depersis} suggest that, if the control input $\bs{u}$ is affected by an exosystem-generated additive disturbance $\bs{d}$, a suitably modified internal-model-based version of the controller   \eqref{eq:new.control} could still guarantee the convergence to the Nash equilibrium in spite of the disturbance. We leave these lines of investigation for future research.
}

\subsection{On the algebraic connectivity}

%\red{... move the previous statements about the strong algebraic connectivity in \cite[Th. 3]{pavel} to here.}

We note that {for the unconstrained case,
\cite[Th. 1]{pavel} establishes 
asymptotic stability without a condition on the graph algebraic connectivity but under a monotonicity assumption on the extended pseudo-gradient $\bs{F}$, an assumption which  is stronger than the Lipschitz continuity considered in Assumption \ref{asspt:lipschitz.extended} (cf., \cite[Remark 1]{pavel}). }
%the strong algebraic connectivity is relaxed via a singular...
{On the other hand, under the same assumption on the Lipschitz continuity of the extended pseudo-gradient $\bs{F}$ considered in our Theorem \ref{convergence1},} \cite[Th. 3]{pavel} relaxes the condition on the algebraic connectivity using a singular perturbation approach. The result, which guarantees exponential stability,  requires however the use of a  gain $\frac{1}{\epsilon}$ which is global to all the agents and must be larger than a bound $\frac{1}{\epsilon^*}$. See \cite[Remark 4]{pavel} for additional discussion on the singular perturbation approach to the problem. The approach we propose in Theorem \ref{convergence1} removes the need of a global parameter. 
{
Similarly, for the constrained case, compared with \cite[Th. 5]{pavel}, 
our Theorem \ref{thm:conv.constr} relaxes the need of a condition on the algebraic connectivity.
%and the discussion following that show that convergence can be achieved via a  controller depending on a global gain parameter $\frac{1}{\epsilon}$ provided that the algebraic connectivity is larger than a quantity  proportional to $\epsilon$. Again,  or the use of a global gain in the controller. 
%Finally, Theorem \ref{thm:conv.constr+coupling} extends the previous results to the case of coupling constraints, showing convergence to generalized Nash equilibria. 
}

\section{Conclusion}
%
%\sgcomment{I would write the conclusion as follows}
%
%{
We have considered the problem to compute a Nash equilibrium in a distributed fashion, with no requirements on the algebraic connectivity of the information exchange graph, nor high-gain assumptions. We have shown that the problem is solvable via dynamic (rather than static) consensus and pseudo-gradient dynamics.
Our analysis and design are based on passivity and an invariance principle for projected dynamical systems. 
Addressing the problem for Nash games with coupling constraints is left as future research.
%}

%We have proposed a new control algorithm that computes Nash equilibria in a distributed fashion. The proposed algorithm improves previous solutions by dispensing with requirements on the algebraic connectivity of the graph used for inter-agent communication and high-gain controllers. Passivity is the  key property that has been exploited to design and analyze the algorithm. The analysis is carried out using an invariance principle for projected dynamical systems. An interesting technical extension of these results concerns the class of games with coupling constrains.  

\appendix

{
\section{Proof of Theorem \ref{th:pavel}} \label{app:proof-th1}
}

For the stability analysis, we consider the Lyapunov function
\[
W(\boldsymbol{x}, \boldsymbol{k})=\tfrac{1}{2} \|\boldsymbol{x}-\overline{\boldsymbol{x}}\|^2 +\tfrac{1}{2} 
\left\| \bs{k} - \overline{\bs{k}} \right\|_{ \Gamma^{-1} }^2
\]
where $\overline{\boldsymbol{x}}=\bs{1}_N \otimes {x}^*$  and $\overline{\boldsymbol{k}} = k^* \bs{1}_N ,
$ with $k^*\in \mathbb{R}$ to determine. 
The time derivative of $W$ is then written as
\be\label{lie.der.W}
\ba{rl}
\dot W(\boldsymbol{x}, \boldsymbol{k}) =& -(\boldsymbol{x}-\overline{\boldsymbol{x}})^\top \mathcal{R}^\top \boldsymbol{F}(\boldsymbol{x}) \\
& - (\boldsymbol{x}-\overline{\boldsymbol{x}})^\top(LKL\otimes I_{M})\boldsymbol{x} \\
& + (\boldsymbol{k} - \overline{\boldsymbol{k}}) ^\top \Gamma^{-1} D(\boldsymbol{\rho})^\top (\Gamma\otimes I_{M}) \bs{\rho}.
\ea\ee

Now, we show that the first addend is bounded as follows:
\begin{multline}
\label{fund.bound}
-(\boldsymbol{x}-\overline{\boldsymbol{x}})^\top \mathcal{R}^\top \boldsymbol{F}(\boldsymbol{x}) \leq \\
-
\left[
\begin{matrix}
\left\| \left( P_N \otimes I_{M} \right) \bs{x} \right\| \\
\left\|  {\rm avg}(\bs{x}) \!-\! \bs{x}^*	 \right\|
\end{matrix}
\right]^\top 
\left[
\begin{matrix}
-\ell_F & \ell \\
\ell & \mu
\end{matrix}
\right]
\left[
\begin{matrix}
\left\| \left( P_N \otimes I_{M} \right) \bs{x} \right\| \\
\left\|  {\rm avg}(\bs{x}) \!-\! \bs{x}^*	 \right\|
\end{matrix}
\right],
\end{multline}

where $P_N := I_N -\tfrac{1}{N} \bs{1}_N\bs{1}_N^\top$ and 
${\rm avg}(\bs{x}):= \tfrac{1}{N} \textstyle \sum_{i=1}^{N} \bs{x}^i = (\tfrac{1}{N} \bs{1}_N^\top \otimes I_{M} )\boldsymbol{x}$.
In fact, following \cite{pavel} with minor modifications, let us define $\boldsymbol{x}^\perp:= (P_N \otimes I_{M}) \boldsymbol{x}$ and $\boldsymbol{x}^\parallel:=  \bs{1}_N \otimes {\rm avg}(\boldsymbol{x})$. Then, the first addend in \eqref{lie.der.W} reads as 
\[\ba{rcl}
- (\boldsymbol{x}-\overline{\boldsymbol{x}})^\top \mathcal{R}^\top \boldsymbol{F}(\boldsymbol{x}) &=& 
- {\boldsymbol{x}^\perp}^\top \mathcal{R}^\top  (\boldsymbol{F}(\boldsymbol{x}) \!-\! \boldsymbol{F}(\boldsymbol{x}^\parallel))\\[2mm]
&& - {\boldsymbol{x}^\perp}^\top \mathcal{R}^\top  \boldsymbol{F}(\boldsymbol{x}^\parallel)\\[2mm]
&& - (\boldsymbol{x}^\parallel-\overline{\boldsymbol{x}})^\top \mathcal{R}^\top  (\boldsymbol{F}(\boldsymbol{x}) \!-\! \boldsymbol{F}(\boldsymbol{x}^\parallel))\\[2mm]
&& - (\boldsymbol{x}^\parallel-\overline{\boldsymbol{x}})^\top \mathcal{R}^\top  \boldsymbol{F}(\boldsymbol{x}^\parallel).
\ea
\]

Since $\| \mathcal{R}^\top\|=1$, by Assumption \ref{asspt:lipschitz.extended}, we have
\be\label{bound1}
\left\| {\boldsymbol{x}^\perp}^\top
\mathcal{R}^\top  (\boldsymbol{F}(\boldsymbol{x}) \!-\! \boldsymbol{F}(\boldsymbol{x}^\parallel))\right\| \le \ell_{\bs{F}} \|(P_N\otimes I_{M}) \boldsymbol{x}\|^2. 
\ee

Then, we note that 
$\boldsymbol{F}(\boldsymbol{x}^\parallel)= F({\rm avg}(\boldsymbol{x}))$. Similarly, since $\overline{\boldsymbol{x}}=\bs{1}_N \otimes x^*$, it holds that $\boldsymbol{F}(\overline{\boldsymbol{x}})= F({x}^*) =\bs{0}$. 
Thus, by the $\ell_F$-Lipschitz continuity of $F$, it holds that 
\be\label{bound2}
\ba{rl}
-{\boldsymbol{x}^\perp}^\top \mathcal{R}^\top  \boldsymbol{F}(\boldsymbol{x}^\parallel)
= & -{\boldsymbol{x}^\perp}^\top \mathcal{R}^\top \left( {F}({\rm avg}(\boldsymbol{x})) - {F}({x}^*) \right) \\
\leq & 
 \ell_{F} \|{\boldsymbol{x}^\perp} \|  \|{\rm avg}(\boldsymbol{x})- {x}^*\|. 
\ea
\ee

Furthermore, we note that 
$\mathcal{R}\boldsymbol{x}^\parallel = {\rm avg}(\boldsymbol{x})$, $\mathcal{R}\overline{\boldsymbol{x}} = {x}^*$,
and therefore, by Assumption \ref{asspt:lipschitz.extended}, the following inequality is true
\begin{equation} \label{bound3}
\begin{array}{l}
- (\boldsymbol{x}^\parallel-\overline{\boldsymbol{x}})^\top \mathcal{R}^\top  (\boldsymbol{F}(\boldsymbol{x}) - \boldsymbol{F}(\boldsymbol{x}^\parallel)) \\
\quad = - ({\rm avg}(\boldsymbol{x})-{x}^*)^\top   (\boldsymbol{F}(\boldsymbol{x}) - \boldsymbol{F}(\boldsymbol{x}^\parallel))  \\
\quad \le \ell_{\bs{F}} \|{\rm avg}(\boldsymbol{x})- {x}^*\| \|{\boldsymbol{x}^\perp}\|.
\end{array}
\end{equation}

By the strong monotonicity stated in Assumption \ref{ass:strong.monotonicity}, we obtain that 
\be\label{bound4}
\ba{l}
- (\boldsymbol{x}^\parallel-\overline{\boldsymbol{x}})^\top \mathcal{R}^\top  \boldsymbol{F}(\boldsymbol{x}^\parallel)
\\
\quad =
 - (\boldsymbol{x}^\parallel-\overline{\boldsymbol{x}})^\top \mathcal{R}^\top  (\boldsymbol{F}(\boldsymbol{x}^\parallel)- 
 \boldsymbol{F}(\overline{\boldsymbol{x}}))\\
\quad =  - ({\rm avg}(\boldsymbol{x}) - {x}^*)^\top  \left( F({\rm avg}(\boldsymbol{x})) - F({x}^*) \right)\\
\quad \le  
 -\mu \|{\rm avg}(\boldsymbol{x}) - {x}^*\|^2.
\ea
\ee

The second addend on the right-hand side in \eqref{lie.der.W} can be rewritten as 
$ -\boldsymbol{x}^\top(LKL\otimes I_{M})\boldsymbol{x}$.
Finally, we rewrite the third addend in \eqref{lie.der.W} as
\[\ba{rcl}
(\boldsymbol{k}{-\overline{\boldsymbol{k}}})^\top \Gamma^{-1} D(\boldsymbol{\rho})^\top (\Gamma\otimes I_{M}) \bs{\rho}
&=& 
\dst \textstyle \sum_{i=1}^N (k_i - k^* ) {\bs{\rho}^i}^\top \bs{\rho}^i\\[4mm]
&=& \bs{\rho}^\top ( (K-K^*) \otimes I_{{M}}) \bs{\rho} \\
&=& \boldsymbol{x}^\top ( L (K-K^*) L \otimes I_{{M}} )\boldsymbol{x},
\ea
\]
where $K^* := k^* I_N$. Thus, the sum of the second and the third addends is equal to 
$-\boldsymbol{x}^\top ( L {K}^* L\otimes I_{M}) \boldsymbol{x}$, from which \cite{kim.depersis}, we conclude that 
\be\label{bound5}
\boldsymbol{x}^\top (L {K}^* L\otimes I_{M})\boldsymbol{x} \ge k^* \lambda_2(L)^2 {\|(P_N \otimes I_{M}) \boldsymbol{x}\|}^2. 
\ee

By replacing the bounds \eqref{bound1}$-$\eqref{bound5} in \eqref{lie.der.W}, we obtain the following inequality:
%\cmargin{changed notation $M$ for matrix to $\mathcal{M}$}
\begin{equation} \label{Lyapunov.in}
\dot W(\boldsymbol{x}, \boldsymbol{k}) \le 
-
\left\| 
\begin{bmatrix}
\|(P_N\otimes I_{M}) \boldsymbol{x}\|
\\
\|{\rm avg}(\boldsymbol{x})- {x}^*\|
\end{bmatrix}
\right\|_{{\mathcal{M}}}^{2}
\end{equation}
where
\begin{equation}
\label{eq:M}
{\mathcal{M}} = 
\begin{bmatrix}
-\ell_{\bs{F}} + k^* \lambda_2(L)^2 & \ell \\
\ell & \mu
\end{bmatrix}
\end{equation}
and the free design parameter $k^*$ is chosen such that ${\mathcal{M}}\succ 0$.

Since the Lyapunov function $W$ is radially unbounded, the inequality in \eqref{Lyapunov.in} shows boundedness of the solutions and convergence to the largest invariant set where $(P_N\otimes I_{M}) \boldsymbol{x} = \bs{0}$ and ${\rm avg}(\boldsymbol{x}) = {x}^*$. On this invariant set, we have $\boldsymbol{x}=\bs{1}_N \otimes {x}^*$. Thus, $ \lim_{t\to \infty} \boldsymbol{x}(t)= \bs{1}_N \otimes {x}^*$. 
\hfill $\blacksquare$

{
\section{Proof of Theorem \ref{invariance.pds}} \label{app:th2}
}

By the compactness of $K\cap \Psi$, the convexity of the set $K$, and the continuous differentiability of $f$, $f$ is Lipschitz continuous on  $K\cap \Psi$. Thus, by the invariance of $K\cap \Psi$ under $\dot x= \Pi_K(f(x))$, for any initial condition $x_0$ in $K\cap \Psi$, there exists a unique Carath\`{e}odory solution defined on $[0, \infty)$ that remains in $K\cap \Psi$. The derivative  $\dot V(x(t))$ exists for almost all $t$ because $x(t)$ is absolutely continuous, and by \eqref{lasalle.ineq}, it satisfies $\dot V(x(t))\le 0$ for almost all $t$. Since $V(x(t))$ is absolutely continuous and $x(t)$ belongs to the compact set $K\cap \Psi$, then  $V(x(t))$ is bounded from below, and the last property, along with $\dot V(x(t))\le 0$ for almost all $t$, implies that $\lim_{t\to \infty} V(x(t; x_0))= V_*$, for some $V_*\in \mathbb{R}$. Consider now a point $p$ in the limit set $\Lambda(x_0)$. Note that $\Lambda(x_0)$ is non-empty because $x(t; x_0)$ is bounded, as a consequence of the Bolzano-Weierstrass theorem. Then, by definition of limit set, $V(p) = V_*$, and since $p$ is a generic point in $\Lambda(x_0)$, $V(p)=V_*$  for all $p\in \Lambda(x_0)$. It is also known that $\lim_{{t\to \infty}} {\rm dist}(x(t;x_0), \Lambda(x_0))=0$. \\
By Lemma \ref{lem:inv.limit.set}, any solution $x(t, p)$ with $p\in \Lambda(x_0)$ remains in $\Lambda(x_0)$ by the Lipschitz continuity of $f$ on $K\cap \Psi$. Thus, since $V(x)$ is constant on $\Lambda(x_0)$, we have $0=\dot V(x(t, p))= \nabla V(x)^\top \Pi_K(f(x(t, p)))$ for almost all $t\in \mathbb{R}_{\ge 0}$. {Now, since} the function $\nabla V(x)^\top \Pi_K(f(x))$ is continuous, {we have}
\begin{multline}
0 = \lim_{t\to 0^+}  \nabla V(x)^\top \Pi_K(f(x(t, p))) \\ 
= \nabla V(x)^\top \Pi_K(f(p)).
\end{multline}
Since $p$ is a generic point in $\Lambda(x_0)$, the equality above shows that $\nabla V(x)^\top \Pi_K(f(p))=0$ for all $p\in \Lambda(x_0)$.  We conclude that $\Lambda(x_0)$ is contained in the largest invariant set  contained in $\{ y \in K\cap \Psi\,\mid\, \nabla V(y)^\top \Pi_K(f(y))=0\}$. 
The proof then follows since $\dst \lim_{t \to 0} {\rm dist}(x(t;x_0), \Lambda(x_0))=0$.
\hfill $\blacksquare$

{
\section{Proof of Theorem \ref{thm:conv.constr}} \label{app:th3}
}

In view of Lemma \ref{lem:diss.in.pds.cl} and Assumption \ref{asspt:lipschitz.extended}, similarly to the proof of Theorem \ref{convergence1}, we obtain that
\begin{multline}
\label{Lyapunov.in.constr2}
\nabla W( \bs{x}, \bs{k} )^\top \Pi_{\bs{\Xi}}\left( \col( \bs{x},\bs{k} ) \,,\, \bs{g}( \bs{x},\bs{k} ) \right) \\ \leq
- \left\| \left[ 
\begin{smallmatrix}
\left\| ( P_N \otimes I_{nN} ) \bs{x}  \right\| \\
\left\| {\rm avg}(\bs{x}) - {x}^* \right\|
\end{smallmatrix}
  \right] \right\|_{\mathcal{M}}^2
\end{multline}
where $\mathcal{M}$ is as in \eqref{eq:M} and $M\succ 0$ for large enough $k^*$.

Having fixed $k^*$, let $\Psi$ be any compact sublevel set of the function $W$ which contains the initial condition in $\bs{\Xi}$. The intersection $\boldsymbol{\Xi}\cap \Psi$ is a compact convex set and therefore $\boldsymbol{g}$ is Lipschitz continuous on it. The last property and the inequality in \eqref{Lyapunov.in.constr2} imply that there exists a unique Carath\`{e}odory solution to \eqref{eq:overall-proj-dynamics-compact-cl}, which belongs to $\boldsymbol{\Xi}\cap \Psi$ for all time and, by Theorem \ref{invariance.pds}, converges to the largest invariant set contained in $\{(\boldsymbol{x}, \bs{k})\in \boldsymbol{\Xi} \cap \Psi: \dot W(\boldsymbol{x}, \bs{k})=0\}$, where, by an abuse of notation, $\dot W(\boldsymbol{x}, \bs{k})$ denotes the right-hand side of \eqref{Lyapunov.in.constr2}. As in Theorem \ref{convergence1}, on this invariant set, we have $\boldsymbol{x}=\bs{1}_N \otimes {x}^*$, which yields the thesis. 
\hfill $\blacksquare$

\bibliographystyle{plain}
\bibliography{int-nask-seeking-2}

\end{document}

%% file: latinte.tex
\newtheorem{theorem}{Theorem}

\newtheorem{definition}{Definition}
\newtheorem{assumption}{Assumption}

\input amssym.def
\input amssym

\def\blacksquare{\hbox{\vrule width 4pt height 4pt depth 0pt}}
\def\square{\hbox{\vrule\vbox{\hrule\phantom{o}\hrule}\vrule}}

\def\col{{\rm col}}
\def\dom{{\rm dom}}
\def\proj{{\rm proj}}

\newcommand{\R}{\mathbb{R}}

\newcommand{\mc}{\mathcal}

\newcommand{\bs}{\boldsymbol}

\DeclareSymbolFont{bbold}{U}{bbold}{m}{n}
\DeclareSymbolFontAlphabet{\mathbbold}{bbold}

\def\v#1{#1}    %** vettore non sottolineato **%

\def \cv {\v  c}

\def\ba#1{\bar #1}

\def\1D#1#2{{{\partial}\over{\partial #2}}#1}

\def\1d#1#2{{{\rm d}\over{{\rm d} #2}}#1}

\newcommand{\dst}{\displaystyle}

\def\be{\begin{equation}}
\def\ee{\end{equation}}

\def\ba{\begin{array}}
\def\ea{\end{array}}

\def\eqa{\begin{eqnarray}}
\def\eqe{\end{eqnarray}}